\documentclass{amsart}
\usepackage{graphicx}
\usepackage{amscd}
\usepackage{amsmath}
\usepackage{amsfonts}
\usepackage{amssymb}
\theoremstyle{plain}

\numberwithin{equation}{section}

\begin{document}
\title[Addendum]{Addendum to Isoperimetry and Symmetrization for Logarithmic Sobolev inequalities}
\author{Joaquim Mart\'{\i}n$^{\ast}$}
\address{Department of Mathematics\\
Universitat Aut\`{o}noma de Barcelona}
\email{jmartin@mat.uab.cat}
\author{Mario Milman}
\address{Department of Mathematics\\
Florida Atlantic University}
\email{extrapol@bellsouth.net}
\urladdr{http://www.math.fau.edu/milman}
\thanks{$^{\ast}$ Supported in part by MTM2004-02299 and by CURE 2005SGR00556}
\thanks{This paper is in final form and no version of it will be submitted for
publication elsewhere.}

\begin{abstract}
We give complete details on an alternative formulation of the
P\'{o}lya-Szeg\"{o} principle that was mentioned in Remark 1 of \cite{mm}. We
also provide an alternative proof \ to a result in \cite{mm}.
\end{abstract}\maketitle

\section{Alternative formulation of P\'{o}lya-Szeg\"{o}}

We discuss in detail the inequality%
\begin{equation}
\int_{0}^{t}((-f_{\mu}^{\ast})^{\prime}(.)I(.))^{\ast}(s)ds\leq\int_{0}%
^{t}\left|  \nabla f\right|  _{\mu}^{\ast}(s)ds,\label{uno}%
\end{equation}
where $g_{\mu}^{\ast}$ denotes the rearrangement of $g$ with respect to the
Gaussian measure $\mu=d\gamma_{n}(x),$ which appears in Remark 1 of our paper
\cite{mm}. In the text we remark that (\ref{uno}) is a reformulation of the
P\'{o}lya-Szeg\"{o} principle and that inequality (4.3) which states that for
any r.i. space $X$ we have%
\begin{equation}
\left\|  (-f^{\ast})^{\prime}(s)I(s)\right\|  _{X}\leq\left\|  \left|  \nabla
f\right|  _{\mu}\right\|  _{X}.\tag{4.3}%
\end{equation}
is a direct consequece of (\ref{uno}). 

We now provide the simple proof.

\begin{proof}
We start with the P\'{o}lya-Szeg\"{o} principle\footnote{This formulation
appears in several places in the literature, in fact, in the Gaussian case it
is implicit in the classical work of Erhard (see \cite{E} page 324).} which we
formulate as
\begin{equation}
\int_{0}^{t}\left|  \nabla f^{\circ}\right|  _{\mu}^{\ast}(s)ds\leq\int
_{0}^{t}\left|  \nabla f\right|  _{\mu}^{\ast}(s)ds.\label{dos}%
\end{equation}
Now, for any positive Young's function $A$ we let $s=\Phi(x_{1})$, and find
\begin{align*}
\int_{0}^{1}A\left(  \left(  -f_{\mu}^{\ast}\right)  {^{\prime}}%
(s)I(s)\right)  ds &  =\int_{\mathbb{R}}A(\left(  -f_{\mu}^{\ast}\right)
{^{\prime}}(\Phi(x_{1}))I(\Phi(x_{1}))\left|  \Phi^{\prime}(x_{1})\right|
dx\\
&  =\int_{\mathbb{R}^{n}}A(\left(  -f_{\mu}^{\ast}\right)  {^{\prime}}%
(\Phi(x_{1}))I(\Phi(x_{1}))d\gamma_{n}(x)\\
&  =\int_{\mathbb{R}^{n}}A(\left|  \nabla f^{\circ}(x)\right|  )d\gamma
_{n}(x).
\end{align*}
where in the last step we have used the fact that
\[
\left(  -f_{\mu}^{\ast}\right)  ^{^{\prime}}(\Phi(x_{1}))I(\Phi(x_{1}%
))=\left(  -f_{\mu}^{\ast}\right)  {^{\prime}}(\Phi(x_{1}))\Phi^{\prime}%
(x_{1})=\left|  \nabla f^{\circ}(x)\right|  .
\]
Since $A$ is increasing, then by Bennett-Sharpley \cite{BS}, exercise 3 pag.
88, we have
\[
\int_{\mathbb{R}^{n}}A(\left|  \nabla f^{\circ}(x)\right|  )d\gamma
_{n}(x)=\int_{0}^{1}A\left(  \left|  \nabla f^{\circ}\right|  _{\mu}^{\ast
}(s)\right)  ds.
\]
Thus,
\[
\int_{0}^{1}A\left(  \left(  -f_{\mu}^{\ast}\right)  {^{\prime}}%
(s)I(s)\right)  ds=\int_{0}^{1}A\left(  \left|  \nabla f^{\circ}\right|
_{\mu}^{\ast}(s)\right)  ds.
\]
Therefore (by \cite{BS}, exercise 5 pag. 88)
\[
\int_{0}^{t}((-f_{\mu}^{\ast})^{\prime}(.)I(.))^{\ast}(s)ds=\int_{0}%
^{t}\left(  \left|  \nabla f\right|  _{\mu}^{\ast}(\cdot)\right)  ^{\ast
}(s)ds,
\]
the second rearrangement is respect to the Lebesgue measure, therefore since
$\left|  \nabla f\right|  _{\mu}^{\ast}(s)$ is decreasing
\[
\left(  \left|  \nabla f\right|  _{\mu}^{\ast}(\cdot)\right)  ^{\ast
}(s)=\left|  \nabla f\right|  _{\mu}^{\ast}(s).
\]
If we combine the previous computation with P\'{o}lya-Szeg\"{o} formulated as
(\ref{dos}) we get
\[
\int_{0}^{t}((-f^{\ast})^{\prime}(.)I(.))^{\ast}(s)ds=\int_{0}^{t}\left|
\nabla f^{\circ}\right|  _{\mu}^{\ast}(s)ds\leq\int_{0}^{t}\left|  \nabla
f\right|  _{\mu}^{\ast}(s)ds
\]
as we wished to show.

In regards to the inequality (4.3) of \cite{mm}: we simply apply the
Hardy-Calder\'{o}n principle to (\ref{uno}) and obtain that for any r.i. space
$X$ we have%
\[
\left\|  (-f^{\ast})^{\prime}(s)I(s)\right\|  _{X}\leq\left\|  \left|  \nabla
f\right|  _{\mu}\right\|  _{X}.
\]
\end{proof}

\bigskip

\section{On the proof of (4.2)}

Professor Andrea Cianchi has kindly brought to our attention that the proof of
(4.2) given in Section 4 of \cite{mm} may not be complete. The problem could
lie in an argument that was not explicitly provided in the text: more
precisely, in the part of the argument when we say that we follow ``Talenti's
argument'' \cite{Ta}. 

In this respect we note that a complete discussion, with proofs, of the
argument in question is given in \cite{ta1}.

Furthermore, although the validity of the proof of (4.2) does not affect the
main results of the paper, we thought it would be prudent to post an
alternative proof independent of this argument. The alternative proof we give
below uses in fact an argument originally given by Professor Cianchi in
\cite{cia}, combined with a suitable twist. We note that we have been familiar
with Cianchi's argument and indeed had occasion to use it in previous
occasions (cf. our recent paper \cite{mm1}, which was cited in \cite{mm}).

We start with the Mazy'a-Talenti inequality: (for a function $h$ we indicate
with $h_{\mu}^{\ast}$ rearrangement with respect to $d\mu=d\gamma_{n}(x))$
\[
\left(  -f_{\mu}^{\ast}\right)  ^{^{\prime}}(s)I(s)\leq\frac{\partial
}{\partial s}\int_{\left\{  \left|  f\right|  >f_{\mu}^{\ast}(s)\right\}
}\left|  \nabla f(x)\right|  d\gamma_{n}(x).
\]
Let us consider a finite family of intervals $\left(  a_{i},b_{i}\right)  ,$
$i=1,\ldots,m$, with $0<a_{1}<b_{1}\leq a_{2}<b_{2}\leq\cdots\leq a_{m}%
<b_{m}<1,$ then
\begin{align*}
\int_{\cup_{1\leq i\leq m}(a_{i},b_{i})}\left(  -f_{\mu}^{\ast}\right)
^{^{\prime}}(s)I(s)ds  &  \leq\int_{\cup_{1\leq i\leq m}(a_{i},b_{i})}\left(
\frac{\partial}{\partial s}\int_{\left\{  \left|  f\right|  >f_{\mu}^{\ast
}(s)\right\}  }\left|  \nabla f(x)\right|  d\gamma_{n}(x)\right)  ds\\
&  =\sum_{i=1}^{m}\int_{\left\{  f_{\mu}^{\ast}(b_{i})<\left|  f\right|  \leq
f_{\mu}^{\ast}(a_{i})\right\}  }\left|  \nabla f(x)\right|  d\gamma_{n}(x)\\
&  =\sum_{i=1}^{m}\int_{\left\{  f_{\mu}^{\ast}(b_{i})<\left|  f\right|
<f_{\mu}^{\ast}(a_{i})\right\}  }\left|  \nabla f(x)\right|  d\gamma_{n}(x)\\
&  =\int_{\cup_{1\leq i\leq m}\left\{  f_{\mu}^{\ast}(b_{i})<\left|  f\right|
<f_{\mu}^{\ast}(a_{i})\right\}  }\left|  \nabla f(x)\right|  d\gamma_{n}(x)\\
&  \leq\int_{0}^{\sum_{i=1}^{m}\left(  b_{i}-a_{i}\right)  }\left|  \nabla
f\right|  _{\mu}^{\ast}(s)ds.
\end{align*}
Now by a routine limiting process we can show that for any measurable set
$E\subset$ $(0,1),$ we have
\[
\int_{E}(-f_{\mu}^{\ast})^{\prime}(s)I(s)ds\leq\int_{0}^{|E|}\left|  \nabla
f\right|  _{\mu}^{\ast}(s)ds.
\]
Therefore
\begin{equation}
\int_{0}^{t}((-f_{\mu}^{\ast})^{\prime}(\cdot)I(\cdot))^{\ast}(s)ds\leq
\int_{0}^{t}\left(  \left|  \nabla f\right|  _{\mu}^{\ast}(\cdot)\right)
^{\ast}(s)ds, \label{aa}%
\end{equation}
where the second rearrangement is respect to the Lebesgue measure. Now, since
$\left|  \nabla f\right|  _{\mu}^{\ast}(s)$ is decreasing, we have
\[
\left(  \left|  \nabla f\right|  _{\mu}^{\ast}(\cdot)\right)  ^{\ast
}(s)=\left|  \nabla f\right|  _{\mu}^{\ast}(s),
\]
and thus (\ref{aa}) yields
\[
\int_{0}^{t}((-f_{\mu}^{\ast})^{\prime}(\cdot)I(\cdot))^{\ast}(s)ds\leq
\int_{0}^{t}\left|  \nabla f\right|  _{\mu}^{\ast}(s)ds,
\]
which in turn is equivalent to the validity of
\begin{equation}
\int_{0}^{1}A\left(  (-f_{\mu}^{\ast})^{\prime}(s)I(s)\right)  ds\leq\int
_{0}^{1}A\left(  \left|  \nabla f\right|  _{\mu}^{\ast}(s)\right)  ds
\label{bb}%
\end{equation}
for every positive Young's function $A.$

Letting $s=\Phi(x_{1})$ in (\ref{bb}), we find
\begin{align*}
\int_{0}^{1}A\left(  \left(  -f_{\mu}^{\ast}\right)  ^{^{\prime}%
}(s)I(s)\right)  ds  &  =\int_{\mathbb{R}}A(\left(  -f_{\mu}^{\ast}\right)
^{^{\prime}}(\Phi(x_{1}))I(\Phi(x_{1}))\left|  \Phi^{\prime}(x_{1})\right|
dx\\
&  =\int_{\mathbb{R}^{n}}A(\left(  -f_{\mu}^{\ast}\right)  ^{^{\prime}}%
(\Phi(x_{1}))I(\Phi(x_{1}))d\gamma_{n}(x)\\
&  =\int_{\mathbb{R}^{n}}A(\left|  \nabla f^{\circ}(x)\right|  )d\gamma
_{n}(x).
\end{align*}
Therefore we get that (\ref{aa}) is equivalent to
\begin{align*}
\int_{\mathbb{R}^{n}}A(\left|  \nabla f^{\circ}(x)\right|  )d\gamma_{n}(x)  &
=\int_{0}^{1}A\left(  (-f_{\mu}^{\ast})^{\prime}(s)I(s)\right)  ds\\
&  \leq\int_{0}^{1}A\left(  \left|  \nabla f\right|  _{\mu}^{\ast}(s)\right)
ds\\
&  =\int_{\mathbb{R}^{n}}A(\left|  \nabla f(x)\right|  )d\gamma_{n}(x).
\end{align*}
This yields (cf. \cite{BS}, exercise 5 pag. 88)
\[
\int_{0}^{t}\left|  \nabla f^{\circ}\right|  _{\mu}^{\ast}(s)ds\leq\int
_{0}^{t}\left|  \nabla f\right|  _{\mu}^{\ast}(s)ds.
\]
as we wished to show.

\end{document}